\documentclass{amsart}  \usepackage{amssymb,latexsym}
\usepackage{graphicx}
\usepackage{etex, verbatim}
\usepackage{color}
\usepackage[pagewise]{lineno}
\usepackage{perpage}

\usepackage[margin=1in]{geometry}
\usepackage{array}
\usepackage[all,2cell,ps]{xy}
\newcommand{\edge}[1]{\ar@{-}[#1]}

\newcommand{\cf}{{\mathcal F}}

\newcommand{\lra}{\longrightarrow}

\newcommand{\ra}{ \rightarrow}

\newcommand{\hra}{\hookrightarrow}

\newcommand{\be}{\begin{equation}}
\newcommand{\ee}{\end{equation}}

\newtheorem{guess}{Theorem}[section]
\newcommand{\bth}{\begin{guess}$\!\!\!${\bf }~}
\newcommand{\eeth}{\end{guess}}
\renewcommand{\bar}{\overline}
\newtheorem{propo}[guess]{Proposition}
\newcommand{\bpropo}{\begin{propo}$\!\!\!${\bf }~}
\newcommand{\epropo}{\end{propo}}

\newtheorem{lema}[guess]{Lemma}
\newcommand{\blem}{\begin{lema}$\!\!\!${\bf }~}
\newcommand{\elem}{\end{lema}}

\newtheorem{defe}[guess]{Definition}
\newcommand{\bdefe}{\begin{defe}$\!\!\!${\bf }~}
\newcommand{\edefe}{\end{defe}}

\newtheorem{coro}[guess]{Corollary}
\newcommand{\bcor}{\begin{coro}$\!\!\!${\bf }~}
\newcommand{\ecor}{\end{coro}}

\newtheorem{rema}[guess]{Remark}
\newcommand{\brem}{\begin{rema}$\!\!\!${\bf }~\rm}
\newcommand{\erem}{\end{rema}}

\newtheorem{exam}[guess]{Example}
\newcommand{\beg}{\begin{exam}$\!\!\!${\bf }~\rm}
\newcommand{\eeg}{\end{exam}}

\newcommand{\tG}{\widetilde{G}}
\newcommand{\tB}{\widetilde{B}}

\newcommand{\tT}{\widetilde{T}
}
\newcommand{\tH}{\widetilde{H}}
\newcommand{\tP}{\widetilde{P}}

\renewcommand{\phi}{\varphi}

\begin{document}
\title{Equivariant Grothendieck ring of a complete symmetric variety
  of minimal rank} \author{V. Uma} \address{Department of Mathematics,
  Indian Institute of Technology-Madras, Chennai, India}
\email{vuma@iitm.ac.in}



\begin{abstract}
  We describe the $G$-equivariant Grothendieck ring of a regular
  compactification $X$ of an adjoint symmetric space $G/H$ of minimal
  rank. This extends the results of Brion and Joshua for the
  equivariant Chow ring of wonderful symmetric varieties of minimal
  rank in \cite{bj} and generalizes the results on the regular
  compactification of an adjoint semisimple group in \cite{u1}.
\end{abstract}  

\maketitle 

\thispagestyle{empty}

\section{Introduction} For a semi-simple group $G$ of adjoint type, we
consider a regular compactification $X$ of a symmetric space $G/H$ of
minimal rank i.e., $\mbox{rk}(G/H)=\mbox{rk}(G)-\mbox{rk}(H)$. The
main examples of symmetric spaces of minimal rank are the groups
$G=G\times G/\mbox{diag}(G)$ and the spaces
${PGL}(2n)/{PSp}(2n)$.

We recall here that a normal complete variety $X$ is called an {\it
  equivariant compactification} of $G/H$ if $X$ contains $G/H$ as an
open subvariety and the action of $G$ on $G/H$ extends to $X$. We
further say that $X$ is a {\it regular compactification} if $X$ is an
equivariant compactification which is regular as a $G$-variety (see
\cite[Section 2]{bdp} and \cite[Section 2.1]{Br1}).

Indeed, the regular compactifications of an adjoint symmetric space
are the complete regular embeddings of a symmetric space or the {\it
  complete symmetric varieties} considered by Bifet, De Concini and
Procesi in \cite{bdp}. In particular, $X$ considered above is a
complete symmetric variety of minimal rank.

The wonderful compactification $X^{wond}$ of a symmetric
space $G/H$ constructed by De Concini and Procesi in \cite{DP1}, is the
unique regular compactification of $G/H$ with a unique closed
orbit. These are called the wonderful symmetric varieties.

Wonderful symmetric varieties have been constructed in characteristic
$0$ by De Concini and Procesi in \cite{DP1} and for an arbitrary
characteristic by De Concini and Springer \cite{DS}.  The geometry and
topology of these varieties have been widely studied by means of their
characterizing wonderful properties (see Section \ref{wondsym}).  For
instance the equivariant cohomology of a wonderful group
compactification has been described by Strickland \cite{str}. Also see
\cite{lp} for the study of the equivariant cohomology ring of a
wonderful symmetric variety.

Let $T\subseteq G$ denote a maximal torus containing a maximal torus
$T_{H}\subseteq H$ with Weyl group $W_{H}$. Among the wonderful
symmetric varieties those of minimal rank (see Definition
\ref{wondsymmin}) have even better geometric properties. For instance,
they have only finitely many $T$-fixed points and finitely many
$T$-invariant curves which can be explicitly described (see
\cite[Lemma 2.1.1]{bj}). This enables the study of the cohomology
theory of these varieties using a precise form of the localization
theorem.

In particular, using the description of the $T$-fixed points and
$T$-invariant curves, the $T$-equivariant Chow ring and hence the
$G$-equivariant Chow ring of wonderful symmetric varieties of minimal
rank have been described by Brion and Joshua in \cite{bj}. The
property of finitely many $T$-fixed points and finitely many
$T$-invariant curves has also been exploited by Tchoudjem \cite{Tc} to
describe the cohomology groups of line bundles on wonderful varieties
of minimal rank.

The complete symmetric varieties have been defined and studied in
\cite{bdp} and \cite{DP2} where their equivariant cohomology ring
structure has been described.

For a regular compactification of a reductive algebraic group (which
includes the adjoint semisimple group) the equivariant Chow ring has
been described in \cite[Section 3]{Br1} and its equivariant
Grothendieck ring has been described in \cite{u1}.

In this article our main aim is to give a description of the
$T$-equivariant and the $G$-equivariant Grothendieck ring
of algebraic equivariant vector bundles on the complete symmetric
variety $X$ of minimal rank.

With this aim in view in Lemma \ref{invariant curves complete} we show
that any complete symmetric variety of minimal rank has only finitely
many $T$-fixed points and $T$-invariant curves and describe them
explicitly. This is an extension of \cite[Lemma 2.1.1]{bj} for the
wonderful symmetric variety of minimal rank.  This is also a
generalization of the corresponding description of the
$T\times T$-fixed points and the $T\times T$-invariant curves for a
regular compactification of a reductive algebraic group in
\cite[p. 160]{Br1} which has been used to describe its
$T\times T$-equivariant and the $G\times G$-equivariant Chow ring in
\cite{Br1} and also to describe its $T\times T$-equivariant and
$G\times G$-equivariant Grothendieck ring in \cite{u1}. 

Let $K_{T}(X)$ (resp. $K_{G}(X)$) denote the $T$-equivariant
(resp. $G$-equivariant) Grothendieck ring of $T$-equivariant
(resp. $G$-equivariant) vector bundles on the complete symmetric
variety $X$. Also $R(T)=K_{T}(pt)$ (resp.  $R(G)=K_{G}(pt)$) denotes
the Grothendieck ring of complex representations of $T$
(resp. $G$). The structure morphism $X\lra \mbox{Spec} ~\mathbb{C}$
induces $R(T)$ (resp.  $R(G)$)-module structure on $K_{T}(X)$ (resp.
$K_{G}(X)$) (see \cite[Section 1.2]{u1}). In particular, when $X$ is
the regular compactification of $G=G\times G/\mbox{diag}(G)$ then we
consider $K_{T\times T}(X)$ (resp. $K_{G\times G}(X)$) as
$R(T)\otimes R(T)$ (resp.  $R(G)\otimes R(G)$)-module. We denote by
$K_{G}(X^{wond})$ the Grothendieck ring of $G$-equivariant vector
bundles on the wonderful compactification $X^{wond}$ of $G/H$.

In Proposition \ref{equivminimalranktorus} and Theorem
\ref{equivminimalrank}, we describe the $K_{T}(X)$ and $K_{G}(X)$
using a precise form of localization theorem \cite[Theorem 1.3]{u1}
for $X$ and the description of the $T$-fixed points and the invariant
curves in Lemma \ref{invariant curves complete}. In particular, we
show that $K_{G}(X)\simeq K_{T}(Y)^{W_{H}}$ where $Y$ denotes the
smooth complete toric variety which is the closure of $T/T_{H}$ in
$X$. These results are the K-theoretic analogue of the description of
the equivariant Chow ring of wonderful symmetric varieties of minimal
rank given in \cite[Theorem 2.2.1]{bj}. These also generalize the
analogous results in \cite[Theorem 2.1, Corollary 2.2, Corollary
2.3]{u1} for a regular compactification of a complex reductive
algebraic group (which includes the case of the regular
compactification of an adjoint semisimple group
$G=(G\times G)/diag(G)$) to all complete symmetric varieties of
minimal rank.

Recall that in \cite[Theorem 2.10]{u1} we also obtain a direct sum
decomposition of $K_{G\times G}(X)$ as a $1\otimes R(G)$-algebra and
describe the multiplication of the graded pieces. In \cite[Corollary
2.11]{u1} we further show the existance of a canonical multifiltration
associated to the direct sum decomposition which in turn gives
$K_{G\times G}(X)$ the structure of a $R(T)\otimes
R(G)$-algebra. Moreover, in \cite[Corollary 2.12]{u1} we gave a
geometric interpretation of the pieces occuring in the direct sum
decomposition.

Our aim in this paper is to extend this kind of description to any
complete symmetric variety of minimal rank.

In Theorem \ref{ds} we give a direct sum decomposition of ${K}_{G}(X)$
as a $1\otimes R(H)$-module and also describe the multiplicative
structure. This extends \cite[Theorem 2.10]{u1} where a similar result
was proved for the regular group compactifications. We then give a
canonical multifiltration for ${K}_{G}(X)$ arising out of the direct
sum decomposition and show that each of the filtered pieces has a
$R(T/T_{H})\otimes R(H)$-module structure (see Corollary \ref{mf} and
Proposition \ref{cansubmodules}). In particular, this shows that
${K}_{G}(X)$ has a structure of a $R(T/T_{H})\otimes R(H)$-algebra
(see Proposition \ref{cansubmodules}). This extends the corresponding
results for regular group compactifications in \cite[Corollary
2.11]{u1}. Further, Corollary \ref{gi} extends \cite[Corollary
2.12]{u1}.

In particular, $R(T/T_{H})=R(T/T_{H})\otimes 1$ can be identified with
the subalgebra of ${K}_{G}(X^{wond})$ generated by the isomorphism
classes of the $G$-linearized line bundles corresponding to the
boundary divisors and $R(H)$ can be identified with $K_{G}(G/H)$ which
is the $G$-equivariant Grothendieck ring of the symmetric space $G/H$
of minimal rank.

The rational equivariant cohomology of complete symmetric varieties
have been described by Bifet, De Concini and Procesi \cite{bdp} in
terms of Stanley-Reisner systems. Our results in Section 3 are an
integral version of their results via K-theory and localization
theorem for the complete symmetric varieties of minimal rank.

Let $\tG$ denote the semi-simple simply connected cover of $G$ and
$\tT$ denote its maximal torus which is a lift of $T$ in $\tG$. Then
$\tG$ and $\tT$ act on $X$ via their canonical projections to $G$ and
$T$ respectively. 

In Section 3.1 we extend Theorem \ref{equivminimalrank} to $K_{\tG}(X)$
and use the results of \cite{mer} to describe the ordinary
Grothendieck ring of $X$.

In Section 3.3 we also give a presentation of $K_{G}(X)$ as an algebra
over $K_{G}(X^{wond})$ which extends \cite[Corollary 4.1]{u2}.

\section{Complete symmetric varieties of minimal
  rank}\label{wondsym}\label{Section 2}

\subsection{Symmetric spaces of minimal rank}

We briefly recall the definition and necessary structure of symmetric
spaces of minimal rank. We refer to \cite[Section 3]{br} and \cite{bj}
for details.  We shall consider algebraic varieties and algebraic
groups over the field of complex numbers.

Let $G$ be a connected reductive algebraic group and let
$$\theta:G\lra G$$ be an involutive automorphism. Let
$$H:=G^{\theta}\subset G$$ denote the subgroup of fixed points. 
The homogeneous space $G/H$ is called a {\em symmetric space}.

It is known that $H$ is reductive and $H^0$ is nontrivial unless
$G$ is a $\theta$-split torus i.e., a torus such that $\theta(g)=g^{-1}$
for every $g\in G$.

Recall that any two maximal $\theta$-split subtori of $G$ are
conjugate in $H^0$, and their common dimension is called the {\it
  rank} of the symmetric space $G/H$ denoted by $\mbox{rk}(G/H)$. Also
every maximal $\theta$-fixed subtorus of $G$ is a maximal torus of
$H$. Any two maximal $\theta$-fixed subtori of $G$ are conjugate in
$H^0$ and their common dimension is the rank of $H$ denoted by
$\mbox{rk}(H)$.

Let $T$ be any $\theta$-stable maximal torus of $G$.

Let
$$T^{\theta}:=\{x\in T~\mid~\theta(x)=x\}$$ and
$$T^{-\theta}:=\{x\in T~\mid~\theta(x)=x^{-1}\}.$$ Then we have \be\label{theta} T=T^{\theta}\cdot
T^{-\theta}~~ \mbox{and} ~~T^{\theta}\cap T^{-\theta} ~~\mbox{is
  finite}.\ee

By (\ref{theta}) it follows that in general we have
\be\label{ineq}\mbox{rk}(G/H)\geq \mbox{rk}(G)-\mbox{rk}(H).\ee If
equality holds in (\ref{ineq}) then we say that the symmetric space
$G/H$ is of {\em minimal rank} i.e., if
$\mbox{rk}(G/H)= \mbox{rk}(G)-\mbox{rk}(H)$. This is equivalent to the
condition that $T^{\theta,0}$ is a maximal $\theta$-fixed subtorus of
$G$ and to the condition that $A:=T^{-\theta,0}$ is a maximal
$\theta$-split subtorus of $G$. This is also equivalent to the
condition that any two $\theta$-stable maximal tori of $G$ are
conjugate in $H^0$ (see \cite[Section 3.1]{br}).

We choose a $\theta$-stable Borel subgroup $B$ and a $\theta$-stable
maximal torus $T$ of $B$ and call $(B,T)$ the standard pair. Also
$\theta$ acts on the Weyl group $W$ and the root system $\Phi$ and
stabilizes $\Phi^{+}$, $\Phi^{-}$ and $\Delta$ which denote
respectively the positive roots, negative roots and the simple roots of
the root system $\Phi$.

We shall assume from now on that $G$ is semi-simple and adjoint so that
$\Delta$ is a basis of $X^*(T)$. We call the symmetric space $G/H$
adjoint as well. We shall also assume that $G/H$ is of minimal rank.
Since $G$ is semi-simple and adjoint, $T_H:=T^\theta$ and
$B_H:=B^{\theta}$ are connected.  Further, $H$ is connected,
semi-simple and adjoint. Thus $T_{H}$ is a maximal torus and $B_{H}$ is
a Borel subgroup of $H$. Moreover, the roots of $(H,T_{H})$ are the
restrictions to $T_{H}$ of the roots of $(G,T)$. The Weyl group $W_{H}$
of $(H,T_{H})$ can be identified with $W^{\theta}$ (see \cite[Lemma
5]{br}).

Recall that the centralizer $C_{G}(A)$ is a Levi subgroup $L$ of a
minimal $\theta$-split parabolic subgroup $P$ of $G$ i.e., $\theta(P)$
is opposite to $P$ and $L=P\cap \theta(P)$.

Let $\Phi_{L}$ denote the root system and $\Delta_{L}\subseteq \Delta$
denote the subset of simple roots of $L$. Under the natural action of
$\theta$ on $\Phi$ we have $\Phi_{L}=\Phi^{\theta}$.  Let
$p:X^*(T)\lra X^*(A)$ denote the restriction map from the character
group of $T$ to the character group of $A$. Then
$p(\Phi)\setminus \{0\}$ is a root system denoted $\Phi_{G/H}$ called
the restricted root system. Also
$\Delta_{G/H}:=p(\Delta\setminus \Delta_{L})$ is a basis of
$\Phi_{G/H}$. The simple restricted roots $\Delta_{G/H}$ can further
be identified with $\alpha-\theta(\alpha)$ for
$\alpha\in \Delta\setminus \Delta_{L}$ under $p$.

The restricted Weyl group $W_{G/H}=N_{G}(A)/C_{G}(A)$ can further
be identified with $N_{W}(A)/C_W(A)$. This yields the exact sequence
$$1\lra W_{L}\lra W^{\theta}=W_H\lra W_{G/H}\lra 1$$  (see
\cite[Subsections 1.2 and 1.3]{bj}).

Let $X^{wond}$ denote the {\em wonderful compactification} of the
adjoint symmetric space $G/H$ constructed by De Concini and Procesi
\cite{DP1} and De Concini and Springer \cite{DS}, which has the
following properties:

\begin{enumerate}
\item $X^{wond}$ is a nonsingular projective variety.

\item $G$ acts on $X^{wond}$ with  an open orbit isomorphic to $G/H$.

\item The complement of the open orbit is the union of
  $r=\mbox{rk}(G/H)$ nonsingular prime divisors $D_1,\ldots, D_r$ with
  normal crossings called the {\em boundary divisors}.

\item The $G$-orbit closures in $X^{wond}$ are exactly the partial
  intersections $\displaystyle X^{wond}_I:=\bigcap_{i\in I} D_i$ where $I$
  runs over the subsets of $\{1,\ldots, r\}$. In particular, each
  $X^{wond}_I$ is smooth.

   \item The unique closed orbit $D_1\cap\cdots\cap D_r$ is isomorphic
     to $G/P\simeq G/\theta(P)$. 
\end{enumerate}

\bdefe\label{wondsymmin} When $G/H$ is of minimal rank, we call the
wonderful compactification $X^{wond}$ a {\em wonderful symmetric variety of
  minimal rank}.\edefe

The wonderful compactification of $G\times G/\mbox{diag}(G)$ is an
example of a wonderful symmetric variety of minimal rank. For other
examples and the complete classification of wonderful symmetric
varieties we refer to \cite[Example 1.4.4]{bj}.

Let $X^{wond}$ be the wonderful symmetric variety of minimal rank.
The associated toric variety $Y^{wond}$ is the closure in $X^{wond}$
of $T/T_{H}$. Then $Y^{wond}$ is a smooth toric variety associated
with the Weyl chambers of the restricted root system $\Phi_{G/H}$. Let
$Y^{wond}_0$ denote the open affine toric subvariety associated with
the positive Weyl chamber dual to $\Delta_{G/H}$. Then
$Y^{wond}=W_{G/H} \cdot Y^{wond}_0$. The distinguished point $z$ of
the unique closed orbit $G/P$ is the unique $T$-fixed point of
$Y^{wond}_0$. Since $W_{G/H}$ acts simply transitively on the
$T$-fixed points, the $T$-fixed points of $Y^{wond}$ are precisely
$w\cdot z$ for $w\in W_{H}/W_{L}\simeq W_{G/H}$. Also the $T$-fixed
points of $X^{wond}$ are precisely $w\cdot z$ for $w\in W/W_{L}$. It
is known that any wonderful variety of minimal rank admits finitely
many $T$-stable curves \cite[Section 10]{Tc}. The precise description
of the $T$-stable curves in the wonderful symmetric variety of minimal
rank $X^{wond}$ as well as that of the $T$-stable curves lying in the
associated toric variety $Y^{wond}$ is given in \cite[Lemma
2.1.1]{bj}. We recall this description below.

\blem\label{invariant curves}\cite[Lemma 2.1.1]{bj}
\begin{enumerate}

\item The $T$-fixed points in $X^{wond}$ (resp., $Y^{wond}$) are exactly the points
  $w\cdot z$ where $w\in W$ (resp., $W_H$), and $z$ denotes the unique
  $T$-fixed point of $Y^{wond}_0$. These fixed points are parametrized by
  $W_G/W_L$ (resp., $W_H/W_L\simeq W_{G/H}$).

\item For any $\alpha\in \Phi^+\setminus \Phi_L^+$, there exists a
  unique irreducible $T$-stable curve $C_{z,\alpha}$ which contains
  $z$ and on which $T$ acts through its character $\alpha$. The
  $T$-fixed points in $C_{z,\alpha}$ are exactly $z$ and
  $s_{\alpha}\cdot z$.

\item For any $\gamma=\alpha-\theta(\alpha)\in \Delta_{G/H}$, there
  exists a unique irreducible $T$-stable curve $C_{z,\gamma}$ on which
  $T$-acts through its character $\gamma$. The $T$-fixed points in
  $C_{z,\gamma}$ are exactly $z$ and $s_{\alpha}s_{\theta(\alpha)}\cdot z$.

\item The irreducible $T$-stable curves in $X^{wond}$ are the
  $W$-translates of the curves $C_{z,\alpha}$ and
  $C_{z,\gamma}$. They are all isomorphic to $\mathbb{P}^1$.

\item The irreducible $T$-stable curves in $Y^{wond}$ are the
  $W_{G/H}$-translates of the curves $C_{z,\gamma}$.
  
\end{enumerate}  
\elem

We now recall the definition of a regular $G$-variety due to Bifet, De
Concini and Procesi \cite{bdp}. (See \cite[Section 1.4]{Br1})

\bdefe\label{regular} A $G$-variety $X$ is said to be regular if it
satisfies the following conditions:

(i) $X$ is smooth and contains a dense $G$-orbit $X^0_{G}$ whose
complement is a union of irreducible smooth divisors with normal
crossings (the boundary divisors).

(ii) Any $G$-orbit closure in $X$ is the transversal intersection of
the boundary divisors which contain it.

(iii) For any $x\in X$, the normal space $T_{x}X/T_{x}(Gx)$ contains a
dense orbit of the isotropy group $G_x$.

\edefe

Consider the adjoint symmetric space $G/H$ of minimal rank. Then $G/H$
is a homogeneous space under $G$ for the standard action from the left
with base point $1\cdot H$. A normal complete variety $X$ is an {\it
  equivariant compactification} of $G/H$ if $X$ contains $G/H$ as an
open subvariety and the action of $G$ on $G/H$ by left multiplication
extends to $X$. We say that $X$ is a {\it regular compactification} of
$G/H$ if $X$ is an equivariant compactification which is regular as a
$G$-variety. The canonical wonderful compactification $X^{wond}$ of
$G/H$ is the unique regular compactification of $G/H$ with a unique
closed orbit.

We have a complete classification of the regular compactifications of
$G/H$ in terms of smooth torus embeddings $Y_0$ of $T/T_{H}$ lying
over $Y_0^{wond}\simeq \mathbb{A}^r$ such that the map
$Y_0\lra Y_0^{wond}$ is proper (\cite[Section 10, Definition
23]{bdp}). Here $Y_0$ is the $T/T_{H}$-toric variety associated to a
fan $\mathcal{F}_+$ which is a smooth subdivision of the positive Weyl
chamber $\mathcal{C}_+$ in the lattice
$X_*(T/T_{H})\otimes \mathbb{R}$ generated by the coweight vectors
dual to the simple restricted roots $\Delta_{G/H}$. Recall that
$W_{G/H}$ acts on the coweight lattice $X_*(T/T_{H})$ by reflection
about the walls of the Weyl chambers. The closure of $T/T_{H}$ in $X$
is the smooth complete toric variety $Y$ which is associated to the
fan $\cf$ in $X_*(T/T_{H})\otimes \mathbb{R}$ whose cones are
$W_{G/H}$-translates of the cones in $\cf_+$.

Let $\cf_{+}(r)$ denote the maximal dimensional cones of $\cf_{+}$
which parametrize the closed $G$-orbits in $X$. For
$\sigma\in \cf_{+}(r)$ we denote by $Z_{\sigma}\simeq G/P$ the
corresponding closed orbit with base point $z_{\sigma}$ (see
\cite[pp. 20-22, Proposition 24]{bdp}).

We have the following precise description of the $T$-fixed points and
the $T$-stable curves in the complete symmetric variety of minimal
rank $X$ as well as that of the $T$-fixed points and the $T$-stable
curves lying in the associated toric variety $Y$.

\blem\label{invariant curves complete}
\begin{enumerate}

\item The $T$-fixed points in $X$ (resp., $Y$) are exactly the points
  $w\cdot z_{\sigma}$ where $w\in W$ (resp., $W_H$), and $z_{\sigma}$
  denotes the $T$-fixed point of $Y_0$ corresponding to the maximal
  cone $\sigma\in \cf_+(r)$. These fixed points are parametrized by
  $\cf_{+}(r)\times W/W_{L}$ (resp., $\cf_{+}(r)\times W_H/W_L$)
  where $W_H/W_L\simeq W_{G/H}$.

\item For any $\alpha\in \Phi^+\setminus \Phi_L^+$, there exists a
  unique irreducible $T$-stable curve $C_{z_{\sigma},\alpha}$ which
  contains $z_{\sigma}$ and on which $T$ acts through its character
  $\alpha$. The $T$-fixed points in $C_{z_{\sigma},\alpha}$ are
  exactly $z_{\sigma}$ and $s_{\alpha}\cdot z_{\sigma}$ for every
  $\sigma\in \cf_{+}(r)$.

\item For any $\gamma=\alpha-\theta(\alpha)\in \Delta_{G/H}$, there
  exists a unique irreducible $T$-stable curve $C_{z_{\sigma},\gamma}$
  on which $T$-acts through its character $\gamma$. The $T$-fixed
  points in $C_{z_{\sigma},\gamma}$ are exactly $z_{\sigma}$ and
  $s_{\alpha}s_{\theta(\alpha)}\cdot z_{\sigma}$. In this case the
  cone $\sigma\in \cf_{+}(r)$ has a facet orthogonal to $\gamma$.

\item There is a unique irreducible $T$-stable curve
  $C_{z_{\sigma}, z_{\sigma'}}$ in $Y_0\subseteq Y\subseteq X$ namely
  the projective line joining $z_{\sigma}$ and $z_{\sigma'}$, which
  are respectively the base points of the distinct closed orbits
  $Z_{\sigma}$ and $Z_{\sigma'}$, whenever the cones $\sigma$ and
  $\sigma'$ in $\cf_{+}(r)$ have a common facet.

\item The irreducible $T$-stable curves in $X$ are the
  $W$-translates of the curves $C_{z_{\sigma},\alpha}$ for every
  $\sigma\in \cf_+(r)$, $C_{z_{\sigma},\gamma}$ whenever
  $\sigma\in \cf_+(r)$ has a facet orthogonal to $\gamma$ and
 $C_{z_{\sigma}, z_{\sigma'}}$ whenever the cones $\sigma$ and
  $\sigma'$ in $\cf_{+}(r)$ have a common facet. They are all
  isomorphic to $\mathbb{P}^1$.

\item The irreducible $T$-stable curves in $Y$ are the
  $W_{G/H}$-translates of the curves $C_{z_{\sigma},\gamma}$ whenever
  $\sigma\in \cf_+(r)$ has a facet orthogonal to $\gamma$ and
  $C_{z_{\sigma}, z_{\sigma'}}$ whenever the cones $\sigma$ and
  $\sigma'$ in $\cf_{+}(r)$ have a common facet.  
\end{enumerate}
\elem

{\bf Proof:} By \cite{bdp} there exists a $G$-equivariant (hence
$T$-equivariant) morphism from $\Phi:X\lra X^{wond}$ which restricts
to the morphism of toric varieties $\phi:Y\lra Y^{wond}$ induced by
the map of fans $\cf_+\lra \mathcal{C}_+$.

Let $x$ be a $T$-fixed point of $X$. Since $\Phi$ is $T$-equivariant
$\Phi(x)$ is a $T$-fixed point of $X^{wond}$. Thus $\Phi(x)=w\cdot z$
for some $w\in W$. Here $z\in Y_0^{wond}$ is the base point of the
unique closed orbit of $X^{wond}$. Since $\Phi$ is $G$-equivariant
$\Phi(w^{-1}\cdot x)=z$. Now $\Phi^{-1}(Y_0^{wond})=Y_0$. Thus
$w^{-1}\cdot x$ is a $T$-fixed point of $Y_0$ and hence is equal to
$z_{\sigma}$ for some $\sigma\in \cf_+(r)$. Thus it follows that
$x=w\cdot z_{\sigma}$ for some $w\in W$ and $\sigma\in
\cf_+(r)$. Since $W_{L}$ acts trivially on $Y$ it follows that the
$T$-fixed points of $X$ are parametrized by $\cf_+(r)\times
W/W_{L}$. This proves (1) of Lemma \ref{invariant curves complete}.

Recall that $W_{H}$ acts on $Y$ and that the $T$-fixed point of $Y$
are $w\cdot z_{\sigma}$ for $\sigma\in\cf_+(r)$ and $w\in W_{H}$ is a
representative of the coset of $W_{H}/W_{L}$. Thus there are
$\cf_+(r)\cdot |W_{H}|/|W_{L}|$ fixed points in $Y$.

Consider the translation $v\cdot Y$ where $v\in W$ is a representative
of the coset $W/W_{H}$. Here $v\cdot Y$ is an irreducible variety
isomorphic to $Y$ with an appropriate twist for the $T$-action.  Then
the $T$-fixed points in $v\cdot Y$ are $v\cdot w\cdot z_{\sigma}$ for
$\sigma\in \cf_{+}(r)$ and $w\in W_{H}/W_{L}$.

Since the $T$-fixed points of $X$ are parametrized by
$\cf_+(r)\times W/W_{L}$ where
$W/W_{L}=\bigcup_{\bar{v}\in W/W_{H}} \bar{v}\cdot W_{H}/W_{L}$. It
follows that the $T$-fixed points of $v\cdot Y$ as $v$ varies over the
coset representatives of $W/W_{H}$ are all distinct and exhaust all
the $T$-fixed points of $X$.

In particular, this implies that $v\cdot Y$ as $v$ varies over the
coset representatives of $W/W_{H}$ are disjoint subvarieties of
$X$. For, if they intersect then the intersection being a complete
$T$-variety will contain a $T$-fixed point which is a contradiction to
the above observation that the collection of $T$-fixed points in
$v\cdot Y$ for distinct coset representatives $v$ of $W/W_{H}$ are
distinct.

Thus a $T$-invariant curve $C$ of $X$ is one the following types.

(i) It lies in some $v\cdot Y$, in which case it is translate by
$v\in W$ to a $T$-invariant curve $C'$ in $Y$. If $C'$ lies in $Y_0$
then $C'$ is of the form (4) in Lemma \ref{invariant curves
  complete}. Thus $C$ is translate by $v$ of a curve of the form (4)
in Lemma \ref{invariant curves complete}. Or else if $C'$ does not lie
in $Y_0$ then $C'$ is conjugate by $w\in W_{H}/W_{L}$ to a
$T$-invariant curve in $Y$ joining $z_{\sigma}$ and
$w\cdot z_{\sigma}$.  Thus $\Phi(C')$ is a $T$-invariant curve of
$Y^{wond}$ joining the two distinct $T$-fixed points $z$ and
$w\cdot z$. Thus $\Phi(C')$ is of the form (3) of Lemma \ref{invariant
  curves}. Since $\Phi$ is $G$-equivariant and hence $T$-equivariant
it follows that $C'$ is of the form (3) of Lemma \ref{invariant curves
  complete}. Hence $C$ is a translate by $w$ of a curve of the form
(3) of Lemma \ref{invariant curves complete}

(ii) It is translate by $v'\in W$ of a $T$-invariant curve $C'$
joining a $T$-fixed point of $Y$ with a $T$-fixed point of $v\cdot Y$
for some coset representative $v$ of $W/W_{H}$ different from
$1$. Since $\Phi(v\cdot Y)=v\cdot Y^{wond}$ for every $v$ the
projection $\Phi(C')$ is a $T$-invariant curve in $X^{wond}$ which
joins $z$ with $v\cdot z$.  Therefore $\Phi(C')$ is of the form (2) of
Lemma \ref{invariant curves}. Since $\Phi$ is $G$-equivariant and
hence $T$-equivariant, $C'$ is of the form (2) in Lemma \ref{invariant
  curves complete}. Thus $C$ is a translate by $v'\in W$ of a curve of
the form (2) in Lemma \ref{invariant curves complete}.  This proves
(5) and (6) of Lemma \ref{invariant curves complete}. $\Box$

We denote by
$\iota_{\sigma}: K_{T}(X)\lra K_{T}(Z_{\sigma})\simeq K_{T}(G/P)$ the
restriction map. For $f\in K_{T}(Z_{\sigma})$ we denote by $f_{w}$ the
restriction of $f$ to the fixed point $w\cdot z_{\sigma}$ for
$w\in W/W_{L}$.

\section{The $G$-equivariant Grothendieck ring of $X$}

Let $G/H$ be a symmetric space of minimal rank. Let $X^{wond}$ denote
the canonical wonderful symmetric variety of minimal rank. Let $X$ be
a projective regular compactification of $G/H$ which in addition is
characterized by the fact that the map of $T/T_{H}$-toric varieties
$Y_0\lra Y_0^{wond}\simeq \mathbb{A}^r$ is a projective morphism.

The following proposition describes the $T$-equivariant Grothendieck
ring of $X$. Also see \cite[Theorem 2.1]{u1} for the corresponding
result for the case of regular group compactifications.

\bpropo\label{equivminimalranktorus} The map
\[\prod_{\sigma\in \cf_{+}(r)} \iota_{\sigma}:K_{T}(X)\lra
  \prod_{\sigma\in \cf_{+}(r)} K_{T}(G/P) \hra R(T)^{|W/W_{L}|\cdot
    |\cf_{+}(r)|}\] is injective and its image consists of all
families $(f_{\sigma})_{\sigma\in \cf_{+}(r)}$ in $K_{T}(G/P)$ such that:

{\em (i)}
$f_{\sigma, w}\equiv f_{\sigma, w\cdot
  s_{\alpha}s_{\theta(\alpha)}}\pmod {1-e^{-w(\gamma)}}$ whenever
$\gamma=\alpha-\theta(\alpha)\in \Delta_{G/H}$ and $w\in W$.

{\em (ii)} $f_{\sigma,w}\equiv f_{\sigma',w}\pmod {1-e^{-w(\chi)}}$ whenever
$\chi\in X^*(T/T_{H})$ and the cones $\sigma$ and $\sigma'$ of
$\cf_{+}(r)$ have a common facet orthogonal to $\chi$ and
$w\in W$.

\epropo {\bf Proof:} By \cite[Theorem 1.3]{u1}, the image of
$K_{T}(X)\lra K_{T}(X^{T})$ is defined by linear congruences
$f_x\equiv f_y\pmod{1-e^{-\chi}}$ whenever $x,y\in X^T$ are connected
by a curve where $T$ acts by the character $\chi$.

Recall that the union of the closed orbits
$\displaystyle\bigcup_{\sigma\in\cf_{+}(r)} Z_{\sigma}$ contains all the $T$-fixed
points of $X$ (see \cite{bdp}, \cite[Lemma 3.4]{br} \cite[Proposition
5.1]{Tc}). Moreover, $T$-fixed points in $Z_{\sigma}$ are the
$W$-translates of $z_{\sigma}$. Since $W_{L}$ fixes $z_{\sigma}$ for
every $\sigma\in \cf_+(r)$ they are parametrized by $W/W_{L}$. In
particular, $z_{\sigma}$ is the distinguished point in $Z_{\sigma}$
which is fixed by $P$. We therefore have the inclusions
$\displaystyle\prod_{\sigma\in \cf_{+}(r)}\iota_{\sigma}:K_{T}(X)\hra
\prod_{\sigma\in \cf_{+}(r)}K_{T}(Z_{\sigma})=K_{T}(G/P)\hra
(R(T))^{|W/W_{L}|\cdot |\cf_+(r)|}$.

By Lemma \ref{invariant curves complete}(5), the invariant curves in
$X$ are the $W$-translates of curves of type (2), (3) or (4). Now, the
$W$-translates of curves of type (2) which are $C_{z_{\sigma},\alpha}$
for $\alpha\in \Phi^{+}\setminus \Phi^{+}_{L}$ lie inside the closed
orbit $Z_{\sigma}$ for $\sigma\in \cf_{+}(r)$ and define its
$T$-equivariant $K$-ring inside $(R(T))^{|W/W_{L}|}$. Thus the image
of $K_{T}(X)$ under $\displaystyle\prod_{\sigma\in \cf_{+}(r)}\iota_{\sigma}$
consists of $\displaystyle (f_{\sigma})\in \prod_{\sigma\in \cf_+(r)}K_{T}(G/P)$
satisfying the congruences (i) and (ii) corresponding respectively to
the $W$-translates of the curves of type (3) and (4). Note that $T$ acts
by $w(\gamma)$ (resp. $w(\chi)$) on the translation by $w\in W$ of the
curve $C_{z_{\sigma},\gamma}$ (resp. $C_{z_{\sigma},
  z_{\sigma'}}$). $\Box$

The $G$-equivariant Chow ring of $X^{wond}$ was described in \cite[Theorem
2.2.1]{bj} using the description of the $T$-fixed points and the
$T$-stable curves given in Lemma \ref{invariant curves} .

We state below the corresponding theorem for the $G$-equivariant
Grothendieck ring $K_{G}(X)$ where $X$ is a regular compactification
of a symmetric space $G/H$ of minimal rank. The proof follows along
similar lines as that of the equivariant Chow ring of $X^{wond}$ by
using \cite[Theorem 1.3]{u1} in place of \cite[Section 3.4]{Br2} and
by the description of the $T$-fixed points and $T$-stable curves of
$X$ given in Lemma \ref{invariant curves complete}. It further
generalizes \cite[Corollary 2.2, Corollary 2.3]{u1} to any complete
symmetric variety of minimal rank.

\bth\label{equivminimalrank} \begin{enumerate} \item[(1)] The ring
  $K_{G}(X)$ consists in all families
  $(f_{\sigma})_{\sigma \in \cf_+(r)}$ of elements in $R(T)^{W_{L}}$
  satisfying \begin{enumerate} \item[$(i)$]
  $s_{\alpha}s_{\theta(\alpha)} \cdot f_{\sigma}\equiv f_{\sigma}\pmod
    {1-e^{-\gamma}}$ whenever
    $\gamma=\alpha-\theta(\alpha)\in \Delta_{G/H}$.\\
  \item[$(ii)$] $f_{\sigma}\equiv f_{\sigma'}\pmod {1-e^{-\chi}}$
    whenever $\chi\in X^*(T/T_{H})$ and the cones $\sigma$ and
    $\sigma'$ of $\cf_{+}(r)$ have a common facet orthogonal to
    $\chi$.\end{enumerate}

\item[(2)] The map
\[r:K_{G}(X)\lra K_{T}(X)^{W}\lra K_{T}(X)^{W_{H}}\lra
  K_{T}(Y)^{W_{H}} \] obtained by composing the canonical maps is an
isomorphism of $R(G)$-algebras where $R(G)=R(T)^{W}$.
\end{enumerate}
\eeth

{\bf Proof:} (1) Note that by \cite[Theorem 1.8]{u1} the restriction
homomorphism $K_{G}(X)\lra K_{T}(X)$ induces an isomorphism
$K_{G}(X)\simeq K_{T}(X)^{W}$. Again by \cite[Theorem 1.8]{u1} the
ring $K_{G}(Z_{\sigma})=K_{G}(G/P)$ is isomorphic to
$K_{T}(G/P)^{W}$. It is further isomorphic to $R(P)=R(T)^{W_{L}}$ via
restriction to $z_{\sigma}$.

Moreover, if $f\in K_{T}(Z_{\sigma})^{W}$, then for each $w\in W$ we
have \be\label{invcong} w\cdot f_{\sigma,1}=f_{\sigma,w}\ee where
$f_{\sigma,1}$ denotes the restriction of $f$ to $z_{\sigma}$ and
$f_{\sigma,w}$ denotes the restriction of $f$ to $w\cdot
z_{\sigma}$. Furthermore, since $W_{L}$ acts trivially on the $T$-fixed
points (\ref{invcong}) implies that $w\cdot f_{\sigma,1}=f_{\sigma,1}$
for $w\in W_{L}$.

Recall that $T$ acts by $w(\gamma)$ (resp. $w(\chi)$) on the
translation of $C_{z,\gamma}$ (resp.  $C_{z_{\sigma}, z_{\sigma'}}$)
by $w\in W$. Thus for
$\displaystyle (f_{\sigma})\in
\prod_{\sigma\in\cf_+(r)}K_{T}(G/P)^{W}$ the congruences (i) and (ii)
can be rewritten as 

(i') $w\cdot f_{\sigma,1}\equiv w\cdot
s_{\alpha}s_{\theta(\alpha)} \cdot f_{\sigma, 1} \pmod {1-e^{-w(\gamma)}}$
whenever $\gamma=\alpha-\theta(\alpha)\in \Delta_{G/H}$ and
$w\in W$.

(ii')
$w\cdot f_{\sigma,1}\equiv w\cdot f_{\sigma',1}\pmod {1-e^{-w(\chi)}}$
whenever $\chi\in X^*(T/T_{H})$ and the cones $\sigma$ and $\sigma'$
of $\cf_{+}(r)$ have a common facet orthogonal to $\chi$ and $w\in W$.

Moreover, since (i') and (ii') are consequences of the congruences
$(i)$ and $(ii)$, this proves (1).

(2) Recall that $T$-fixed points of $Y$ are the $W_{H}$ translates of
$z_{\sigma}$ for $\sigma\in \cf_{+}(r)$. Also $W_L$ acts trivially on
the $T$-fixed points.  Furthermore, by Lemma \ref{invariant curves
  complete} (6) and by observing that $T$ acts by $w(\gamma)$
(resp. $w(\chi)$) on the translation of the curve
$C_{z_{\sigma},\gamma}$ (resp.  $C_{z_{\sigma}, z_{\sigma'}}$) by
$w\in W_{H}$, it can be seen that $K_{T}(Y)$ can be identified with
the the tuples
$\displaystyle (f_{\sigma,w})_{w\in W_{H}/W_{L}, \sigma\in
  \cf_{+}(r)}$ in $R(T)^{|W_{H}/W_{L}|\cdot |\cf_{+}(r)|}$ satisfying
the congruences

(i'')
$f_{\sigma, w}\equiv f_{\sigma, w \cdot
  s_{\alpha}s_{\theta(\alpha)}}\pmod {1-e^{-w(\gamma)}}$ whenever
$\gamma=\alpha-\theta(\alpha)\in \Delta_{G/H}$ and $w\in W_{H}$.

(ii'') $f_{\sigma,w}\equiv f_{\sigma',w}\pmod {1-e^{-w(\chi)}}$ whenever
$\chi\in X^*(T/T_{H})$, the cones $\sigma$ and $\sigma'$ of
$\cf_{+}(r)$ have a common facet orthogonal to $\chi$ and
$w\in W_{H}$.

Here $f_{\sigma,w}$ denotes the restriction of $f\in K_{T}(Y)$ to
$w\cdot z_{\sigma}$ for $\sigma\in\cf_{+}(r)$ and $w\in
W_{H}/W_{L}$. Furthermore, if $f\in K_{T}(Y)^{W_{H}}$ then
$f_{\sigma,w}=w\cdot f_{\sigma}$ where $f_{\sigma}$ denotes the
restriction of $f$ to $z_{\sigma}$. Thus it follows that the ring
$K_{T}(Y)^{W_{H}}$ is identified with the tuples
$\displaystyle ( f_{\sigma})_{\sigma\in \cf_{+}(r)}\in
\prod_{\sigma\in \cf_{+}(r)} R(T)^{W_L}$ satisfying the congruences
$(i)$ and $(ii)$. Hence the theorem. $\Box$

We shall denote by $\gamma$ the restricted root
$\alpha-\theta(\alpha)$ in $\Delta_{G/H}$ for
$\alpha\in \Delta_{H}\setminus \Delta_L$ and by $s_{\gamma}$ the
corresponding simple reflection which is the image in $W_{G/H}$ of the
element $s_{\alpha}\cdot s_{\theta(\alpha)}\in W_{H}$. (Recall that
$s_{\gamma}^2=1$ \cite[Definition 6, Lemma 10.2]{Tc}.)

We isolate below the description of $K_{G}(X^{wond})$ as a particular
case of the above theorem. This is analogous to \cite[Lemma 3.2]{u1}
in the case of the wonderful compactification of an adjoint
semi-simple group.

\bcor\label{precise} The ring
$K_{G}(X^{wond})$ is identified with the subring of
$R(T)^{W_{L}}$ defined by the congruences
\[f\equiv s_{\alpha}s_{\theta(\alpha)}\cdot f \pmod {1-e^{\gamma}} \]
for $\gamma=\alpha-\theta(\alpha)\in \Delta_{G/H}$,
$\alpha\in \Delta_{H}\setminus \Delta_{L}$.  \ecor {\bf Proof:} In
this case $|\cf_{+}(r)|=1$ since $X^{wond}$ has a unique closed
$G$-orbit. Thus by Theorem \ref{equivminimalrank} (1)
$K_{G}(X^{wond})$ can be identified with the subring of $R(T)^{W_{L}}$
satisfying the congruences $(i)$. Hence the corollary. $\Box$

Let $\tG$ denote the simply connected cover of $G$ with projection
$\pi:\tG\ra G$ and $\tT:=\pi^{-1}(T)$. Then $\tG$ is a semi-simple
simply connected algebraic group with $\tT$ as a maximal torus. Let
$\tP:=\pi^{-1}(P)$. There exists an involution $\widetilde{\theta}$ of
$\tG$ which induces the involution $\theta$ of the adjoint quotient
$G$. Let $\tH$ denote the fixed points of $\tG$ under
$\widetilde{\theta}$. Then $\tT_{H}:=\tH\cap \tT$ is a maximal torus
of $\tH$.

We have the following theorem analogous to Theorem
\ref{equivminimalrank} obtained by replacing $G$ by $\tG$ and $T$ by
$\tT$ and by considering the action of $\tG$ and $\tT$ on $X$ and $Y$
through their canonical surjections to $G$ and $T$ respectively. We
omit the proof to avoid repetition.

\bth\label{equivminimalranksc}
The map 
\[r:K_{\tG}(X)\lra K_{\tT}(X)^{W}\lra K_{\tT}(X)^{W_{H}}\lra K_{\tT}(Y)^{W_{H}} \]
obtained by composing the canonical maps is an isomorphism.
\eeth

We further have the following corollary analogous to Corollary
\ref{precise}.  \bcor\label{precisesc} The ring $K_{\tG}(X^{wond})$ is
identified with the subring of $R(\tT)^{W_{L}}$ defined by the
congruences
\[f\equiv s_{\alpha}s_{\theta(\alpha)}\cdot f \pmod {1-e^{\gamma}} \]
for $\gamma=\alpha-\theta(\alpha)\in \Delta_{G/H}$,
$\alpha\in \Delta_{H}\setminus \Delta_{L}$.  \ecor

\brem\label{boundary} For $X^{wond}$ the set of simple restricted
roots $\Delta_{G/H}$ are in bijection with the set of boundary
divisors. Let $D_{\gamma}$ denote the boundary divisor corresponding
to $\gamma\in \Delta_{G/H}$.  Let $\mathcal{L}_{\gamma}$ denote the
$G$-linearized (and hence $\tG$-linearized) line bundle on $X$
corresponding to the $G$-stable boundary divisor $D_{\gamma}$ for
$\gamma\in \Delta_{G/H}$. Moreover, $\mathcal{L}_{\gamma}$ has a
$\tG$-invariant section whose zero locus is the boundary divisor
$D_{\gamma}$ for $\gamma\in \Delta_{G/H}$.  Note that $X^*(T/T_{H})$
has a basis consisting of the simple roots $\gamma\in
\Delta_{G/H}$. Thus $R(T/T_{H})=\mathbb{Z}[X^*(T/T_{H})]$ is generated
as a $\mathbb{Z}$-algebra by $e^{\gamma}, \gamma\in
\Delta_{G/H}$. Thus we can identify $R(T/T_{H})$ with the subalgebra
of $K_{G}(X^{wond})$ generated by $[\mathcal{L}_{\gamma}]$ for
$\gamma\in \Delta_{G/H}$ (see \cite[Remark 3.7]{u1}). In particular,
$\mathcal{L}_{\gamma}\mid_{Y^{wond}_0}$ is the $T/T_{H}$-linearized
line bundle on the toric variety $Y^{wond}_0$ corresponding to the ray
in $\cf_+$ generated by $\omega_{\gamma}$.  \erem

\brem\label{extcobord}({\it Extension to other oriented cohomology
  theories}) The methods of \cite{bj} have earlier been used by
V. Kiritchenko and A. Krishna in \cite{KK} to compute the rational
equivariant algebraic cobordism ring of wonderful symmetric varieties
of minimal rank. We wish to remark here that Lemma \ref{invariant
  curves complete} above allows us to extend the description in
\cite[Theorem 6.4]{KK} to all regular symmetric varieties of minimal
rank.  \erem

\subsection{Ordinary K-ring of $X$}

Note that we have the equalities
$R(\tT)^{W_{H}}=K_{\tG}(\tG/\tB)^{W_{H}}$ is an
$R(\tG)=R(\tT)^{W}$-subalgebra of $K_{\tG}(\tG/\tB)=R(\tT)$.

Consider the augmentation map $\epsilon: R(\tG)\lra \mathbb{Z}$  which
takes an element $[V]$ to $\dim(V)$. Thus $\mathbb{Z}$ becomes a
$R(\tG)$-module via the augmentation map.

We have the isomorphisms (\cite[Theorem 6.1.22]{mer})
\be\label{ordmin}K(X)\simeq \mathbb{Z}\otimes_{R(\tG)} K_{\tG}(X)\ee
and
\be\label{ordtor}K_{\tT}(Y)^{W_{H}}=K_{\tG}(\tG\times_{\tB}Y)^{W_H}
.\ee

\brem Here $\tG\times_{\tB} Y$ is a toric bundle with fibre the
$T/T_{H}$-toric variety $Y$ and base $\tG/\tB$.
\erem

\bpropo\label{kord} We have the following isomorphism of ordinary
$K$-rings:
\[K(X)\simeq K(\tG\times_{\tB} Y)^{W_{H}}.\]
\epropo {\bf Proof:} Follows immediately from Theorem
\ref{equivminimalrank}, (\ref{ordmin}) and (\ref{ordtor}).  $\Box$

\subsection{Determination of the structure of $K_{G}(X)$}

We first fix some notations similar to \cite[Section 2.1]{u1}. Let
$\mathcal{F}$ denote the (smooth projective) fan associated to
$Y$. Recall that the Weyl group $W_{H}$ acts on $\mathcal{F}$ by
reflection across the walls of the Weyl chambers and the cones in
$\mathcal{F}$ get permuted by this action of $W_{H}$ and each cone is
stabilized by the reflections corresponding to the walls of the Weyl
chambers on which it lies. Let $W_{\tau}$ denote the subgroup of
$W_{H}$ which fixes the cone $\tau\in \mathcal{F}$. Then, in
particular, $W_{\sigma}=W_{L}$ for $\sigma\in \mathcal{F}(r)$ where
$r=\mbox{rk}(G/H)$ and $W_{\{0\}}=W_{H}$. In particular,
$W_{H}\supseteq W_{\tau}\supseteq W_{L}$ for every $\tau\in
\cf$. Indeed, $W_{H}$ acts on $\cf$ via its projection to
$W_{H}/W_{L}=W_{G/H}$ which is the Weyl group of the restricted root
system with simple roots $\Delta_{G/H}$.

Let $\{\rho_j\mid j=1,\ldots, d\}$ denote the set of edges of the fan
$\mathcal{F}$ and let $\tau(1)$ denote the set of edges of the cone
$\tau$ for every $\tau\in \mathcal{F}$. Let $v_j$ denote the primitive
vector along the edge $\rho_j$. Let $O_{\tau}$ denote the
$T/T_{H}$-orbit in $Y$ corresponding to $\tau\in \mathcal{F}$. Let
$L_j$ denote the $T/T_{H}$-equivariant line bundle on $Y$
corresponding to the edge $\rho_j$. Note that $L_j$ has a
$T/T_{H}$-invariant section $s_j$ whose zero locus is
$\overline{O_{\rho_j}}$. Recall that $Y_0$ is the toric variety
associated to the fan $\mathcal{F}_{+}$ consisting of the cones
$\tau\in \mathcal{F}$ which lie in the positive Weyl chamber.

Let $X_{F}:=\prod_{\rho_j\in F}(1-X_j)$ in the Laurent polynomial
algebra $\mathbb{Z}[X_1^{\pm 1},\ldots, X_d^{\pm1}]$ for every
$F\subseteq \{\rho_j\mid j=1,\ldots, d\}$. In particular we let
$X_\tau:=X_{\tau(1)}=\prod_{\rho_j\in \tau(1)}(1-X_j)$ for every
$\tau\in \mathcal{F}$.

Let
$C_{\tau}:=X_{\tau}\cdot\mathbb{Z}[X_{j}^{\pm1}:\rho_j\in \tau(1)]$.
Recall from \cite[Lemma 2.8]{u1} that we have the additive
decomposition \be\label{dsdec} K_{T/T_{H}}(Y_0)=\bigoplus_{\tau\in
  \cf_{+}} C_{\tau}\ee which follows from the Stanley-Reisner
presentation of the $T/T_{H}$-equivariant Grothendieck ring of $Y_0$,
$K_{T/T_{H}}(Y_0)=\mathbb{Z}[X_j^{\pm 1}:\rho_j\in \cf_{+}(1)]/\langle
X_{F} ~~\forall~~ F\notin \cf_{+}\rangle$ (see \cite[Theorem
6.4]{VV}).

Similarly by \cite[Theorem 6.4]{VV} we also have the following Stanley Reisner
presentation for the $T/T_{H}$-equivariant Grothendieck ring of $Y$:
\be\label{srp} K_{T/T_{H}}(Y)\simeq \mathbb{Z}[X_1^{\pm1},\ldots,
X_d^{\pm1}]/\langle X_{F} ~~\forall~~ F\notin \cf\rangle\ee

Furthermore, since $W_{H}$ acts on $\mathcal{F}$ we have an action of
$W_{H}$ on $\mathbb{Z}[X_1^{\pm 1},\ldots, X_d^{\pm 1}]$, given by
$w(X^{\pm 1}_{\rho_j})=X^{\pm 1}_{w(\rho_j)}$ for every $w\in
W_{H}$. Thus $w(X_{F})=X_{w(F)}$ for $F\subseteq \{\rho_j\mid
j=1,\ldots, d\}$ and $w\in W_{H}$, and since $W_{H}$ permutes the
cones of $\mathcal{F}$ we further get an action of $W_{H}$ on the
Stanley-Reisner algebra $\mathbb{Z}[X_1^{\pm1},\ldots, X_d^{\pm1}]/\langle
X_{F} ~~\forall~~ F\notin \cf\rangle$ so that (\ref{srp})  is an
isomorphism of $W_{H}$-modules where the $W_{H}$-action on
$K_{T/T_{H}}(Y)$ is induced from the $W_{H}$-action on $Y$.

We further have the additive decomposition \[\mathbb{Z}[X_1^{\pm1},\ldots, X_d^{\pm1}]/\langle
X_{F} ~~\forall~~ F\notin \cf\rangle=\bigoplus_{\tau\in \mathcal{F}}X_{\tau}\cdot
\mathbb{Z}[X_j^{\pm1}:\rho_j\in \tau(1)],\]

where $W_{H}$ acts on the right hand side as follows:
\[w\cdot (X_{\tau}\cdot
\mathbb{Z}[X_j^{\pm1}:\rho_j\in \tau(1)])=X_{w(\tau)}\cdot
\mathbb{Z}[X_j^{\pm1}:\rho_j\in w(\tau)(1)] \] for every $w\in W_{H}$.

Thus we can write

\[K_{T/T_{H}}(Y)= \bigoplus_{\tau\in \cf_{+}}\bigoplus_{w\in
    W_{H}/W_{\tau}} X_{w(\tau)}\cdot \mathbb{Z}[X_j^{\pm1}:\rho_j\in w(\tau)(1)]\]
where we recall that $W_{\tau}$ denotes the subgroup of $W_H$ which
fixes $\tau\in\cf_{+}$.

This implies that as a $W_{H}$-module we can write \be\label{ind1}
K_{T/T_{H}}(Y)=\bigoplus Ind_{W_{\tau}}^{W_{H}} C_{\tau}\ee where
$C_{\tau}:=X_{\tau}\cdot\mathbb{Z}[X_{j}^{\pm1}:\rho_j\in
\tau(1)]$. Further, since $C_{\tau}$ is fixed by $W_{\tau}$,
\be\label{ind2}
Ind_{W_{\tau}}^{W_{H}}C_{\tau}=\mathbb{Z}[W_{H}/W_{\tau}]\otimes
C_{\tau}.\ee

Recall (see \cite{mer}) that the structure morphism $Y\lra pt$ induces
canonical inclusions $R(T)\hra K_{T}(Y)$ (resp.
$R(T/T_{H})\hra K_{T/T_{H}}(Y)$) which gives an $R(T)$-algebra (resp.
$R(T/T_{H})$-algebra) structure on $K_{T}(Y)$
(resp. $K_{T/T_{H}}(Y)$).

Note that we have an exact sequence \be\label{exact} 0\lra
X^*(T/T_{H})\lra X^*(T)\lra X^*(T_{H})\lra 0\ee of character groups,
where $X^*(T/T_{H})$ denotes the characters of $T$ that are trivial
when restricted to $T_{H}$. The exact sequence (\ref{exact}) splits by
choosing a basis of $X^*(T_{H})$ and by lifting every element of this
basis to a character of $T$. This implies that we have the following
isomorphism \be\label{splitting}R(T)\simeq R(T/T_H)\otimes R(T_{H})\ee
where $R(T)=\mathbb{Z}[X^*(T)]$, $R(T_{H})=\mathbb{Z}[X^*(T_{H})]$,
$R(T/T_{H})=\mathbb{Z}[X^*(T/T_{H})]$ (see \cite[page 485]{bj} for the
analogous statement for Chow ring).

More generally, we have the following proposition whose proof is
similar to that of \cite[Lemma 1.7]{u1} (see \cite[Lemma 2.3.2]{bj}
for the analogous statement for Chow ring).

\bpropo\label{kiso} We have the following isomorphisms of rings
\begin{enumerate}
\item[(i)] $K_{T}(Y)\simeq K_{T/T_{H}}(Y)\otimes_{\mathbb{Z}} R(T_{H})$ \item[(ii)]
  $K_{T}(Y)\simeq K_{T/T_{H}}(Y)\otimes _{R(T/T_{H})}
  R(T)$ \end{enumerate} \epropo {\bf Proof:} It suffices to prove (i). For,
(ii) will then follow from (\ref{splitting}) since the right hand side of
(ii) can be written as
$K_{T/T_{H}}(Y)\otimes _{R(T/T_{H})\otimes 1} R(T/T_{H})\otimes
R(T_{H})$ which is isomorphic to the right hand side of (i).

We now show (i). Since $\chi\in X^*(T_{H})$ can be lifted to a
character $\chi'\in X^*(T)$ (by (\ref{exact})), we have a canonical
homomorphism of rings $K_{T/T_{H}}(Y)\otimes R(T_{H})\lra K_{T}(Y)$
where $e^{\chi}\in R(T_{H})=\mathbb{Z}[X^*(T_{H})]$ maps to
$e^{\chi'}\in R(T)=\mathbb{Z}[X^*(T)]\hra K_{T}(Y)$ and the map from
$K_{T/T_{H}}(Y)$ to $K_{T}(Y)$ is induced by the surjection
$T\lra T/T_{H}$.

To define the inverse of the above homomorphism, we let $E$ be a
$T$-equivariant vector bundle on $Y$. Since $Y$ is a
$T/T_{H}$-variety, the $T_{H}$-action on $Y$ is trivial. Thus we get a
$T_{H}$-action on every fibre of $E$ which gives a weight space
decomposition on each fibre. Further, the $T$-equivariant vector
bundle $E$ is locally trivial so the weights of the restricted
$T_{H}$-action are locally constant. Moreover, since $Y$ is
irreducible these weights are constant over $Y$. Thus we get an
isotypical decomposition $\displaystyle E=\bigoplus E_i$ where $E_i$
denotes the subbundle of $E$ whose fibre is the eigenspace for the
$T_{H}$-action corresponding to a character $\chi_i$ of $T_{H}$. We
can write $E_i=e^{{\chi'}_i}\otimes E^i$ where $E^i$ is a
$T$-equivariant vector bundle with a trivial $T_{H}$-action and hence
a $T/T_{H}$-equivariant vector bundle on $Y$. We therefore define the
inverse map which sends $[E]\in K_{T}(Y)$ to
$\sum_{i}e^{\chi_i}\otimes [E^i]\in R(T_{H}) \otimes_{\mathbb{Z}}
K_{T/T_{H}}(Y)$.  $\Box$

\brem Note that the above proposition holds in a more general setting
where we can take \[1\lra T'\lra T\lra T''\lra 1\] to be an exact
sequence of tori and $Y$ to be any irreducible $T''$-variety. The
proof follows verbatim. We have stated it in the above form for our
convenience. \erem

  With the above notations we have the following theorem which extends
   \cite[Theorem 2.10]{u1}.

   \bth\label{ds} Let $X$ be a complete symmetric variety of minimal
   rank and $Y$ be the corresponding torus embedding of $T/T_H$.  The
   ring $K_{G}(X)$ has the following direct sum decomposition as a
   $1\otimes R(H)$-module:

\be\label{kdec} K_{G}(X)=\bigoplus_{\tau\in \cf_{+}} C_{\tau}\otimes
R(T_{H})^{W_{\tau}},\ee where $R(H)=R(T_{H})^{W_{H}}$ acts naturally on
  the second factor in each piece of the above decomposition. Further,
  the multiplictive structure of $K_{G}(X)$ can be described from the
  above decomposition as follows: Let $a_{\tau}\otimes b_{\tau}\in
  C_{\tau}\otimes R(T_{H})^{W_{\tau}}$ and let $a_{\sigma}\otimes b_{\sigma}\in
  C_{\sigma}\otimes R(T_{H})^{W_{\sigma}}$. Then

  \[ (a_{\tau}\otimes b_{\tau})\cdot (a_{\sigma}\otimes
    b_{\sigma})=\left\{\begin{array}{ll} a_{\tau}\cdot a_{\sigma}\otimes
                         b_{\tau}\cdot 
                         b_{\sigma}&\mbox{if}~\tau~\mbox{and}~\sigma~\mbox{span~
                                   ~  the~ cone}~\gamma\\ 0&
                                                             \mbox{if}~~\tau~\mbox{and}~\sigma~\mbox{do~not ~span~
                                     a~ cone~in}~\cf_{+}\end{array}\right.  \]

Here $a_{\tau}\cdot a_{\sigma}\otimes b_{\tau}\cdot b_{\sigma}\in
C_{\gamma}\otimes R(T_{H})^{W_{\gamma}}$ when $\tau$ and $\sigma$ span
$\gamma$. The multiplication in the first factor  is as in
$K_{T/T_{H}}(Y_0)$ where $C_{\tau}\cdot C_{\sigma}\subseteq C_{\gamma}$. 
  
\eeth {\bf Proof:} We have the following isomorphism by Theorem
\ref{equivminimalrank}: \[K_{G}(X)\simeq K_{T}(Y)^{W_{H}}.\] This together with
Proposition \ref{kiso}(i) will imply that \be\label{iso1}
K_{G}(X)\simeq (K_{T/T_{H}}(Y)\otimes_{\mathbb{Z}}
R(T_{H}))^{W_{H}}.\ee

Now, by (\ref{ind1}) and (\ref{ind2}) we get:
\[K_{T/T_{H}}(Y)\otimes_{\mathbb{Z}}
R(T_{H})=\bigoplus_{\tau\in \cf_{+}} \mathbb{Z}[W_{H}/W_{\tau}]\otimes
C_{\tau}\otimes R(T_{H}).\]

Further, by taking $W_{H}$-invariants on either side we get:

\be\label{iso2} (K_{T/T_{H}}(Y)\otimes_{\mathbb{Z}} R(T_{H}))^{W_{H}}=
  \bigoplus_{\tau\in \cf_{+}} C_{\tau}\otimes R(T_{H})^{W_{\tau}}. \ee

  Thus from (\ref{iso1}) and (\ref{iso2}), the additive decomposition
  (\ref{kdec}) follows.

  Observe that $C_{\tau}\cdot C_{\sigma}\subseteq C_{\gamma}$ whenever
  $\tau$ and $\sigma$ span a cone $\gamma$ in $\cf_{+}$ and
  $C_{\tau}\cdot C_{\sigma}=0$ if $\tau$ and $\sigma$ do not span a
  cone in $\cf_{+}$. Furthermore, whenever
  $\gamma=\langle \tau,\sigma\rangle$ in $\cf_{+}$ the product
  $R(T_{H})^{W_{\tau}}\cdot R(T_{H})^{W_{\sigma}}\subseteq
  R(T)^{W_{\gamma}}$ since $R(T_{H})^{W_{\tau}}$ and
  $R(T_{H})^{W_{\sigma}}$ are both subrings of $R(T)^{W_{\gamma}}$.

  The multiplicative structure now follows exactly as in the proof of
  \cite[Theorem 2.10]{u1}.  $\Box$

  We have the following corollary which extends \cite[Corollary]{u1}.

  \bcor\label{mf} The ring
  $K_{G}(X)\simeq \bigoplus_{\tau\in \cf_{+}} C_{\tau}\otimes
  R(T_{H})^{W_{\tau}}$ admits a multifiltration
  $\{F_{\tau}\}_{\tau\in \cf_{+}}$ where the filtered pieces
  are
  \[F_{\tau}=\bigoplus_{\tau\preceq \sigma}C_{\sigma}\otimes
    R(T_{H})^{W_{\sigma}},\] where $F_{\tau}\supseteq F_{\sigma}$ whenever
  $\tau\preceq \sigma$, and $F_{\{0\}}=K_{G}(X)$. Further, under the
  multiplication described in Theorem \ref{ds}, we have
  $F_{\tau}\cdot F_{\sigma}\subseteq F_{\gamma}$ whenever
  $\gamma=\langle \tau,\sigma \rangle$ in $\cf_{+}$. In particular
  $F_{\{0\}}\cdot F_{\tau}\subseteq F_{\tau}$ for all
  $\tau\in \cf_{+}$.  \ecor {\bf Proof:} The existance of the
  filtration $\{F_{\tau}\}_{\tau\in \cf_{+}}$ and the properties
  follow by definition. Further, since the filtered pieces multiply by
  the multiplication rule given in Theorem \ref{ds}, it follows that
  $F_{\tau}\cdot F_{\sigma}\subseteq F_{\gamma}$ whenever
  $\gamma=\langle \tau,\sigma \rangle$ in $\cf_{+}$ and
  $F_{\tau}\cdot F_{\sigma}\subseteq 0$ whenever $\tau$ and $\sigma$
  do not span a cone in $\cf_{+}$.
$\Box$

In the following proposition we show that the above multifiltration
establishes the existence of the structure of
$K_{T/T_{H}}(Y_0)\otimes R(H)$-algebra on $K_{G}(X)$. We further show
that the multifiltration also implies the existence of certain
canonical $K_{T/T_{H}}(Y_0)\otimes R(H)$-submodules which will be
subsequently used in the next section. As a consequence of Theorem
\ref{ds}, we further show that $K_{G}(X)$ has the structure of a
canonical $K_{T/T_{H}}(Y_0)\otimes R(H)$-submodule of
$K_{T/T_{H}}(Y_0)\otimes R(T_{H})^{W_{L}}$. 

\bpropo\label{cansubmodules} \begin{enumerate}
\item[(i)] We have a canonical inclusion
  $K_{T/T_{H}}(Y_0)\otimes R(H)\subseteq
  F_{\{0\}}=K_{G}(X)$ as $1\otimes R(H)$-subalgebra.

\item[(ii)] We have a canonical inclusion of the
  $K_{T/T_{H}}(Y_0)\otimes R(H)$-submodule
  \[\prod_{\rho_{j}\in \tau(1)}X_{\tau}\cdot K_{T/T_{H}}(Y_0)\otimes
  R(T_{H})^{W_{\tau}}\] in $F_{\tau}\subseteq F_{\{0\}}=K_{G}(X)$.

\item[(iii)] We have a canonical inclusion of $K_{G}(X)$ in
  $K_{T/T_{H}}(Y_0)\otimes R(T_{H})^{W_{L}}$ as a
  $K_{T/T_{H}}(Y_0)\otimes R(H)$-subalgebra. (This is analogous to
  \cite[Proposition 2.5]{u1} and \cite[Proposition 2.1]{u2}.)

\end{enumerate}
    
\epropo {\bf Proof:} (i) From (\ref{dsdec}) it follows that
\be\label{dsdectensor}K_{T/T_{H}}(Y_0)\otimes R(H)=\bigoplus_{\tau\in
  \cf_+}C_{\tau}\otimes R(H).\ee

Also \be\label{inclusion1}C_{\tau}\otimes R(H)=C_{\tau}\otimes R(T_{H})^{W_{H}}\subseteq
C_{\tau}\otimes R(T_{H})^{W_{\tau}}\ee for every $\tau\in \cf_{+}(1)$.

Since both (\ref{kdec}) and (\ref{dsdectensor}) are decopositions as
as $1\otimes R(H)$-algebras, (i) follows from (\ref{inclusion1}).

(ii) Again from (\ref{dsdec}) we have
\be\label{dstensor1}X_{\tau}\cdot K_{T/T_{H}}(Y_0)\otimes R(T_{H})^{W_{\tau}}=\bigoplus_{\sigma\in
  \cf_{+}}X_{\tau}\cdot C_{\sigma}\otimes R(T_{H})^{W_{\tau}}.\ee

Moreover, by the relations in the Stanley Reisner presentation of
$K_{T/T_{H}}(Y_0)$ we have
$X_{\tau}\cdot C_{\sigma}\subseteq C_{\gamma}$ whenever $\tau$ and
$\sigma$ span a cone $\gamma$ in $\cf_+$ and
$X_{\tau}\cdot C_{\sigma}=0$ whenever $\tau$ and $\sigma$ do not span
a cone in $\cf_+$. Thus
\[X_{\tau}\cdot C_{\sigma}\otimes R(T_{H})^{W_{\tau}}\subseteq
  C_{\gamma}\otimes R(T_{H})^{W_{\gamma}}\] since
$R(T)^{W_{\tau}}\subseteq R(T_{H})^{W_{\gamma}}$ whenever
$\tau\preceq \gamma$.

Now, (ii) follows by (\ref{dstensor1}) and by the definition of
$F_{\tau}$ from Corollary \ref{mf}.

(iii) Recall  that for every $\tau\in \cf_+$ we have $W_{\tau}\supseteq
W_{L}$ which further implies that $R(T_{H})^{W_{\tau}}\subseteq
R(T_{H})^{W_{L}}$.

Thus for every $\tau\in \cf_+$ we have the inclusion of
$1\otimes R(H)$-submodules
\[C_{\tau}\otimes R(T_{H})^{W_{\tau}}\subseteq C_{\tau}\otimes R(T_{H})^{W_{L}} \]

Taking direct sum over all $\tau\in \cf_+$ we get:

\be\label{kdecsub} \bigoplus_{\tau\in \cf_+}C_{\tau}\otimes
  R(T_{H})^{W_{\tau}}\subseteq \bigoplus_{\tau\in \cf_+} C_{\tau}\otimes R(T_{H})^{W_{L}} .\ee

We see by (\ref{dsdec}) that the right hand side of the above equation
is isomorphic to $K_{T/T_{H}}(Y_0)\otimes R(T_{H})^{W_{L}}$.

Now, by (\ref{kdec}) and (\ref{kdecsub}) (iii) follows.  $\Box$

In the following corollary we give a geometric interpretation of
Theorem \ref{ds}. This is analogous to \cite[Corollary 2.12]{u1} for
regular compactifications. 

\bcor\label{gi} Let $N_{\tau}\simeq \oplus_{\rho_j\in \tau(1)}L_j$
be the normal bundle of $V_{\tau}=\bar{O_{\tau}}$ in $Y$. Let
$N_{\tau}\mid_{O_{\tau}}$ denote the restriction of the normal bundle
to $O_{\tau}$ so that so that

\[\lambda_{-1}(N_{\tau}\mid_{O_{\tau}}):=\prod_{\rho_j\in
    \tau(1)}(1-[L_j]\mid_{O_{\tau}})\in K_{T/T_{H}}(O_{\tau}).\]
Then we have the following decomposition:

\[K_{G}(X)\simeq \bigoplus_{\tau\in
    \cf_{+}}\lambda_{-1}(N_{\tau}\mid_{O_{\tau}})\cdot
  K_{T/T_{H}}(O_{\tau})\otimes  R(T_{H})^{W_{\tau}}.\]

Let $P_{\tau}:=\lambda_{-1}(N_{\tau}\mid_{O_{\tau}})\cdot
K_{T/T_{H}}(O_{\tau})$ for each $\tau\in \cf_{+}$. Then the above
decoposition is a ring isomorphism where the multiplication on the
right hand side is given as follows:

\[P_\tau\cdot P_{\sigma}\subseteq \left\{\begin{array}{ll}
P_{\gamma}&\mbox{if} ~\tau ~\mbox{and} ~\sigma ~\mbox{span ~the~
            cone}~\gamma~\in ~\cf_{+},\\0 &\mbox{if} ~\tau ~\mbox{and}~
                                         \sigma ~\mbox{do ~not~ span~ a~
            cone~in} ~\cf_{+}

                                         \end{array}\right.\]  

\ecor

The structure of rational equivariant cohomology of regular embeddings
have been described in \cite{bdp} and \cite[Theorem 2.3]{lp}. This
result for a complete symmetric variety of minimal rank can
alternately be obtained by proving results analogous to Theorem
\ref{equivminimalrank} and Theorem \ref{ds} for the equivariant
cohomology of a complete symmetric variety (see \cite{lp} for the
corresponding results on wonderful symmetric varieties).

Let $X_{F}=\prod_{\rho_{j}\in F}X_j$ for every
$F\subseteq \{\rho_j\mid j=1,\ldots,d\}$ in the polynomial ring
$\mathbb{Q}[X_1,\ldots, X_d]$. In particular
$X_{\tau}:=X_{\tau(1)}=\prod_{\rho_j\in \tau(1)}X_{j}$ for every
$\tau\in \cf$. Let $e(N_{\tau})$ denote the equivariant Euler class of
the normal bundle $V_{\tau}=\bar{O_{\tau}}$. Let $S:=H^*_{T_{H}}(pt)$
By \cite[Theorem 8, p.7]{bdp} we have:
\[H^*_{T/T_{H}}(Y)\simeq \mathbb{Q}[X_1,\ldots, X_d]/\langle
  X_{F}~\forall ~F~\notin \cf \rangle\] Let $e(N_{\tau})$ denote the
equivariant Euler class of the normal bundle of
$V(\tau)=\bar{O_{\tau}}$ which is equal to the top Chern class of
$\displaystyle\bigoplus_{\rho_j\in \tau(1)}L_j$.  We then have the
following description of the equivariant cohomology of a complete
symmetric variety of minimal rank (see \cite{lp} and \cite{bdp}). We
consider cohomology ring with $\mathbb{Q}$-coefficients.

\bth\label{equivcohomminimalrank}
\[H^*_{G}(X)\simeq \bigoplus _{\tau\in \cf_+}e(N_{\tau}\mid
  _{O_{\tau}})\cdot H_{T/T_{H}}(O_{\tau})\otimes S^{W_{\tau}}.\]
\eeth

Note that $S^{W_{H}}=H_{H}^*(pt)\simeq H_{G}^*(G/H)$ and
$H_{T/T_{H}}^*(Y_0)$ can be identified with the subalgebra of
$H^*_{G}(X)$ generated by the equivariant classes of the boundary
divisors.  \brem From the above decomposition, an analogue of
Proposition \ref{cansubmodules} for the equivariant cohomology can
also be proved which will show that $H^*_{G}(X)$ is a
$H^*_{T/T_{H}}(Y_0)\otimes S^{W_H}$-submodule of
$H^*_{T/T_{H}}(Y_0)\otimes S^{W_{L}}$ (see \cite[Proposition
2.3.3]{bj} for a similar statement on $A^*_{G}(X)_{\mathbb{Q}}$). Our
aim is to show that $H^*_{G}(X)$ is a free module of rank $|W_{G/H}|$
over $H^*_{T/T_{H}}(Y_0)\otimes S^{W_H}$ and to find a basis of this
free module. This will be an extension of the results of Strickland
\cite{str} for the compactification of adjoint semisimple group to any
complete symmetric variety of minimal rank. This will be taken up in a
future work. We also note here that similar description of the
equivariant cohomology ring of any regular compactification of an
adjoint semisimple group has been obtained by Strickland (see
\cite[Theorem 6.1]{St}). \erem

\subsection{Comparison with wonderful symmetric variety}

Let $X$ be a regular compactification of the adjoint symmetric space
$G/H$ and $X^{wond}$ be the canonical wonderful compactification.

Let $\gamma_1,\ldots, \gamma_r$ denote the simple restricted roots.

Recall that $\mathcal{L}_{\gamma_i}$ are $\tG$-linearized line bundles
on $X^{wond}$ such that $\tP$ operates on $\mathcal{L}_{\gamma_i}$ by
the character $\gamma_i$ for $1\leq  i\leq r$ where $z$ is the base
point of the unique closed orbit.

Furthermore, since the centre $Z$ of $\tG$ acts trivially on
$X^{wond}$ and hence on the fibre by the character $\gamma_i$ the line
bundle $\mathcal{L}_{\gamma_i}$ is actually
$G=\tG/Z$-linearized. Moreover, $\mathcal{L}_{\gamma_i}$ admits a
$G$-invariant section $s_i$ whose zero locus is the boundary divisor
$D_{\gamma_i}$ for $1\leq i\leq r$.

We recall the following construction from \cite[Section 3.5]{bdp} and
\cite[Section 10]{DP2}.

The bundle
$\mathcal{V}=\displaystyle\bigoplus_{1\leq i\leq r}
\mathcal{L}_{\gamma_i}$ being a direct sum of line bundles on
$X^{wond}$ admits a natural action of the $r$-dimensional torus
$\mathbb{G}_m^r$. Put

\[ \mathcal{P}:=\mathcal{V}\setminus\bigcup_{i=1}^r
  \mathcal{L}_{\gamma_1}\oplus \cdot \oplus
  \hat{\mathcal{L}}_{\gamma_i}\oplus \cdots \oplus \mathcal{L}_{\gamma_r}.\]

Then $\mathcal{P}$ is the principal $T/T_{H}=\mathbb{G}_m^r$-bundle
associated with $\mathcal{V}$ over $X^{wond}$. The section
$s=s_1\oplus \cdots \oplus s_r$ of $\mathcal{V}$ clearly maps $G$ to
$\mathcal{P}$. Therefore using $s$ we can embed $X^{wond}$ in the
bundle $\mathcal{P}\times_{T/T_{H}}\mathbb{A}^r$ which is nothing but
$\mathcal{V}$. Since the bundles $\mathcal{L}_{\gamma_i}$ are
$G$-linearized it follows that $\mathcal{P}$ is a left $G$-space
whose bundle map $\mathcal{P}\lra X^{wond}$ is $G$-equivariant
for the canonical $G$-action on $X^{wond}$. Moreover, the right
$T/T_{H}$-action is compatible with the left $G$-action on
$\mathcal{P}$.

Consider the proper morphism of toric varieties
$\phi: Y_0\lra Y_0^{wond}\simeq \mathbb{A}^r$ corresponding to a smooth
subdivision $\cf_+$ of the positive Weyl chamber
$\mathcal{C}_{+}$. Since $Y_0$ is a $T/T_{H}$-toric variety we can
form the associated toric bundle $\mathcal{P}\times _{T/T_{H}} Y_0$
over $X^{wond}$. The morphism $\phi$ further induces
$\bar{\phi}: \mathcal{P}\times_{T/T_{H}} Y_0\lra \mathcal{P}\times
_{T/T_{H}} Y_0^{wond}=\mathcal{V}$ which is isomorphic to the pull
back of $\mathcal{P}\times _{T/T_{H}} Y_0 $ to $\mathcal{V}$ via the
bundle projection $p:\mathcal{V}\lra X^{wond}$.

From \cite[Definition 23]{bdp} and \cite[Section 3.2]{DP2},
$X=X_{Y_0}=\bar{\phi}^{-1}(s(X_{wond}))$ which can be identified with
pull-back of
$\mathcal{P}\times_{T/T_{H}} Y_0\stackrel{\bar{\phi}}{\lra}
\mathcal{P}\times _{T/T_{H}} \mathbb{A}^r$ via the section $s$. Thus
$X=X_{Y_0}=s^*(p^*(\mathcal{P}\times _{T/T_{H}} Y_0))$. Now,
$p\circ s=Id_{X^{wond}}$ hence $s^*\circ p^*=Id^*$. Thus
\be\label{isomorphism} X\simeq \mathcal{P}\times _{T/T_{H}} Y_0.\ee

For $u\in X^*(T/T_{H})$ let
$\mathcal{L}_u:=\mathcal{P}\times_{T/T_{H}} \mathbb{C}_{u}$ denote the
$G$-equivariant line bundle on $X^{wond}$ associated to the character
$u$. Let $[\mathcal{L}_u]_{G}$ denote its class in $K_{G}(X^{wond})$.

Let $X$ be a projective regular compactification of $G/H$ so that
$Y_0$ is a semi-projective $T/T_{H}$-toric variety. Then the following
theorem follows by \cite[Theorem 5.1]{u2} and (\ref{isomorphism}).

\bth\label{relkwond} The $G$-equivariant Grothendieck ring of $X$ has
the following description as an algebra over the $G$-equivariant
Grothendiek ring of $X^{wond}$:

\be \label{kisomorphism} K_{G}(X)\simeq \frac{K_{G}(X^{wond})[X_j^{\pm
    1}~: ~\rho_j\in \cf_{+}(1)]}{\mathcal{J}}\ee where $\mathcal{J}$
is the ideal in $K_{G}(X^{wond})[X_j^{\pm 1}~: ~\rho_j\in \cf_{+}(1)]$
generated by the elements $X_{F}$ for $F\notin \cf_+$ and
$\displaystyle\prod_{\rho_j\in \cf_{+}(1)}{X_j}^{\langle u,
  v_j\rangle}-[\mathcal{L}_u]_{G}$ for $u\in X^*(T/T_{H})$.  \eeth

Recall that the adjoint symmetric spaces of minimal rank are exactly
the products of the symmetric spaces listed in \cite[Examples
1.4.4]{bj}.  In particular, \cite[Examples 1.4.4]{bj}(1) is the
adjoint semisimple group $G=G\times G/\mbox{diag}(G)$ and
\cite[Examples 1.4.4]{bj}(2) is the space $G/H=PGL(2n)/PSp(2n)$. The
$G\times G$-equivariant $K$-ring of (1) has already beed described in
\cite{u1}. In the following example we shall describe the
$G$-equivariant $K$-ring of the wonderful compactification $X^{wond}$
of $PGL(2n)/PSp(2n)$.

\beg\label{Example:}

Consider the group $G=PGL(2n)$ and the involution $\theta$ 
  associated with the symmetry of the Dynkin diagram so that
  $H=PSp(2n)$ and $\mbox{rk}(G/H)=n-1$. We shall consider the case
  when $n=3$ for simplicity. The general case is similar.

Let $\Delta_{G}=\{\alpha_1,\alpha_2,\alpha_3,\alpha_4,\alpha_5\}$
  be the simple roots of the root system $A_{5}$ corresponding to
  $PGL(6)$. Then $\Delta_{L}=\{\alpha_1,\alpha_3,\alpha_5\}$. Also
  $\Delta_{G}\setminus \Delta_{L}=\{\alpha_2,\alpha_4\}$. Thus
  $\Delta_{G/H}=\{\gamma_1=\alpha_2-\theta(\alpha_2),
  \gamma_2=\alpha_4-\theta(\alpha_4)\}$ are the simple restricted
  roots.  

  Since $T_{H}$ is a maximal torus of the adjoint semisimple group of
  type $C_3$ it is known that $R(T_{H})=\mathbb{Z}[X^*(T_{H})]$ where
  $X^*(T_{H})$ has a basis consisting of the $\mbox{rk}(H)=3$ simple
  roots $\Delta_{H}=\{\beta_1=2e_1,\beta_2=e_2-e_1,\beta_3=e_3-e_2\}$
  where $\beta_1$ is a long root and $\beta_2$ and $\beta_3$ are short
  roots.

  Note that $\theta(\alpha_1)=\alpha_1$, $\theta(\alpha_3)=\alpha_3$,
  $\theta(\alpha_5)=\alpha_5$,
  $\theta(\alpha_2)=-\alpha_1-\alpha_2-\alpha_3$ and
  $\theta(\alpha_4)=-\alpha_3-\alpha_4-\alpha_5$.
  
  Thus under the restriction map $q:X^*(T)\lra X^*(T_{H})$ (see
  \cite[Lemma 1.4.2]{bj}) $q^{-1}(\beta_1)$ consists of the unique
  root $\alpha_1$, $q^{-1}(\beta_2)$ consists of the two strongly
  orthogonal roots $\alpha_2+\alpha_3$ and $-\alpha_1-\alpha_2$ and
  $q^{-1}(\beta_3)$ consists of the two strongly orthogonal roots
  $-\alpha_3-\alpha_4$ and $\alpha_4+\alpha_5$.  The restricted simple
  roots are
  $\gamma_1=\alpha_2-\theta(\alpha_2)=\alpha_1+2\alpha_2+\alpha_3$ and
  $\gamma_2=\alpha_4-\theta(\alpha_4)=\alpha_3+2\alpha_4+\alpha_5$.

  Here $W=S_{6}$, $W_{H}$ is the semidirect product of
  $S_2\times S_2\times S_2$ with $S_{3}$,
  $W_{L}=S_2\times S_2\times S_2$ and $W_{G/H}=S_{3}$.

  In this case the toric variety $Y^{wond}$ is associated to the Weyl
  chambers of type $A_{2}$ and $Y^{wond}_0$ is the toric variety associated
  to the positive Weyl chamber in the coweight lattice spanned by the
  fundamental coweights $\omega^{\vee}_{\gamma_1}$ and
  $\omega^{\vee}_{\gamma_2}$ dual to the simple roots $\gamma_1$ and
  $\gamma_2$ respectively.

  The simple reflections in $W_{G/H}$ corresponding to $\gamma_1$ and
  $\gamma_2$ are $s_{\gamma_1}$ and $s_{\gamma_2}$ respectively.  Here
  $\cf_{+}$ consists of $4$ cones namely $\{0\}$, $\rho_1,\rho_2$ and
  $\sigma=\langle \rho_1,\rho_2\rangle$, where the coweight vectors
  $\omega^{\vee}_{\gamma_1}$ and $\omega^{\vee}_{\gamma_2}$ are the
  primitive vectors along the rays $\rho_1$ and $\rho_{2}$
  respectively.  Note that $W_{\sigma}=W_{L}=S_2\times S_2\times S_2$,
  $W_{\{0\}}=W_{H}$ which is the semidirect product of $W_{L}$ with
  $S_{3}$, $W_{\rho_1}$ is the semidirect product of $W_{L}$ with the
  subgroup of $W_{G/H}=S_3$ generated by $s_{\gamma_2}$, $W_{\rho_2}$
  is the semidirect product of $W_{L}$ with the subgroup of
  $W_{G/H}=S_3$ generated by $s_{\gamma_1}$.

  Let $\mathcal{R}:=R(T_{H})^{W_{L}}$. Further,
  $R(PSp(6))=R(H)=R(T_{H})^{W_{H}}=\mathcal{R}^{W_{G/H}}$. Furthermore,
  $K_{T/T_{H}}(Y_0)\simeq R(T/T_{H})=\mathbb{Z}[e^{\pm\gamma_1},e^{\pm
    \gamma_2}]\simeq \mathbb{Z}[X_{\rho_1}^{\pm 1}, X_{\rho_{2}}^{\pm
    1}]$ where $X_{\rho_1}$ (respectively $X_{\rho_2}$)maps to the
  class of the trivial line bundle on $Y_0$ where $T/T_{H}$ acts on
  the fibre via the character $\gamma_1$ (respectively $\gamma_2$)
  . The following is an additive decomposition of $K_{PGL_6}(X)$ as a
  $1\otimes R(PSp(6))$-submodule of $R(T/T_{H})\otimes \mathcal{R}$.
  \begin{align*}&\mathbb{Z}\otimes R(PSp(6))\oplus
    (1-e^{\gamma_1})\cdot
                  \mathbb{Z}[e^{\pm \gamma_1}]\otimes \mathcal{R}^{s_{\gamma_2}}\oplus\\
                &\oplus (1-e^{\gamma_2})\cdot\mathbb{Z}[e^{\pm
                  \gamma_2}]\otimes \mathcal{R}^{s_{\gamma_1}}\oplus
                  (1-e^{\gamma_1})\cdot (1-e^{\gamma_2})\cdot
                  \mathbb{Z}[e^{\pm \gamma_1} , e^{\pm
                  \gamma_2}]\otimes \mathcal{R}\end{align*}

              \eeg

 \brem By the parametrization of line bundles on
              spherical varieties (see \cite[Section 2.2]{br2}) the
              $\tG$-linearized line bundles on $X$ correspond to
              $PL(\cf_+)$ which denotes the piecewise linear functions
              on the fan $\cf_+$. Let $\mathcal{L}_h$ be the
              $\tG$-linearized line bundle on $X$ corresponding to
              $h=(h_{\sigma})_{\sigma\in \cf_+(r)}$ then
              $\tP=\pi^{-1}(P)$ acts on
              $\mathcal{L}\mid_{z_{\sigma}}$ by the character
              $h_{\sigma}$. Thus
              $\mathbf{L}_h=\mathcal{L}_{h}\mid_{Y_0}$ is a
              $\tT/\tT_{H}$-linearized line bundle on $Y_0$
              corresponding to the piecewise linear function
              $h\in PL(\cf_+)$ so that $\tT/\tT_{H}$ acts on
              $\mathbf{L}_h\mid_{z_{\sigma}}$ by the character
              $h_{\sigma}$. Since it is known that
              $\mbox{Pic}^{\tT/\tT_{H}}(Y_0)$ generates
              $K_{T/T_{H}}(Y_0)$ as a ring (see \cite{VV}) it follows
              that $K_{\tT/\tT_{H}}(Y_0)$ is the subring of
              $K_{\tG}(X)$ generated by $\mbox{Pic}^{\tG}(X)$. We can
              further identify $K_{\tG}(G/H)\simeq R(\tH)$. We wonder
              whether $K_{\tG}(X)$ has a canonical
              $K_{\tT/\tT_{H}}(Y_0)\otimes R(\tH)$-module structure so
              that $K_{\tG}(X)$ is a free module of rank $|W_{G/H}|$
              over $K_{\tT/\tT_{H}}(Y_0)\otimes R(\tH)$ and if there
              exists a basis for this free module generalizing the
              results in \cite[Section 3]{u1} and \cite{u2} to the
              setting of all complete symmetric varieties of minimal
              rank.  See \cite[Theorem 4 (ii)]{br} for a similar
              statement describing $K(X^{wond})_{\mathbb{Q}}$ as a
              free module of rank $|W_{H}/W_{L}|$ over the classes of
              the boundary divisors $[\mathcal{O}_{D_{\gamma}}]$ for
              $\gamma\in \Delta_{G/H}$ together with the classes of
              equivariant vector bundles $\Phi(\lambda)$ induced by
              the fundamental weights $\lambda$ of $H$.

              An analogue of \cite[Proposition 2.3.1]{bj} for the
              $\tT$ and $\tG$-equivariant Grothendieck ring may give
              us a way to extend \cite[Theorem 3.3]{u1} to all
              wonderful symmetric varieties of minimal rank. However
              the proof of \cite[Proposition 2.3.1]{bj} does not
              extend to the integral setting of the Grothendieck ring,
              due to the reason that there is no $W_{L}$-invariant
              splitting for short exact sequence of character groups
              \be\label{tildeexact} 0\lra X^*(\tT/\tT_{H})\lra
              X^*(\tT)\lra X^*(\tT_{H})\lra 0.\ee This
              can be seen from the following example of
              $PGL_4/PSp(4)$. \erem

\beg\label{counter} In the case of $G=PGL(4)$ and $H=PSp(4)$,
$\tG=SL(4)$ and $\tH=Sp(4)$.  The maximal torus $\tT$ is the diagonal
torus $(t_1,t_2,t_3,t_4)\in (\mathbb{C}^*)^4$ satisfying
$t_1\cdot t_2\cdot t_3\cdot t_4=1$. The subtorus $\tT_{H}$ of $\tT$
consisting of $(t_1,t_2,t_3,t_4)$ such that $t_1\cdot t_2=1$ and
$t_3\cdot t_4=1$.  Here $W_{L}\simeq S_2\times S_2$ is the subgroup of
$S_4$ generated by the transpositions $\{(1,2) , (3,4)\}$.  The
existance of a $W_{L}$-invariant splitting for the restriction
homomorphism \be\label{res} q:X^*(\tT)\lra X^*(\tT_{H})\ee would mean
the existance of a complementary subtorus of $\tT_{H}$ in $\tT$ which
is $W_{L}$-invariant. A $1$-dimensional subtorus of $\tT$ is a
$1$-parameter subgroup $(t^{a_1},t^{a_2}, t^{a_3}, t^{a_4})$ for
$t\in\mathbb{C}^*$ where $(a_1,a_2,a_3,a_4)\in \mathbb{Z}^4$ satisfies
$a_1+a_2+a_3+a_4=0$. If this is $W_{L}$ invariant then it will be of
the form $(t^{a_1},t^{a_1}, t^{a_3},t^{a_3})$ such that
$a_1+a_3=0$. Thus it will be of the form $(t,t, t^{-1},t^{-1})$ for
$t\in \mathbb{C}^*$. However, this subtorus is not a complement of
$\tT_{H}$ since it intersects $\tT_{H}$ along the subgroup $\{\pm I\}$
of order $2$. Thus in general there cannot exist a $W_{L}$-invariant
splitting for (\ref{res}).  \eeg

\vspace{0.2cm}

\brem({\it Assumptions on the base field}) Although we work over the field of
complex numbers all the results in this paper should hold
over an algebraically closed field $k$ of characteristic $\neq~2$ (see
\cite[p. 473]{bj} and \cite[p. 273-274]{DS}).\erem 

\vspace{0.2cm}

\noindent {\bf Acknowledgments:} I am very grateful to Prof. Michel
Brion for his valuable comments and suggestions while preparing this
manuscript. I thank the referee for a careful reading of the
manuscript and for very helpful comments and suggestions. I wish to
mention that Remark \ref{extcobord} has been added based on the
referee's observation.

\vspace{0.5cm}
\noindent{\bf Statements and Declarations} I hereby declare that this work
has no related financial or non-financial conflict of interests.

\end{document}